*Research Article*

# θ-Metric Space: A Generalization

**Farshid Khojasteh,[1] Erdal Karapınar,[2] and Stojan Radenović[3]**

[1] *Department of Mathematics, Arak Branch, Islamic Azad University, Arak, Iran*
[2] *Department of Mathematics, Atilim University, İncek, 06836 Ankara, Turkey*
[3] *Faculty of Mechanical Engineering, University of Belgrade, Kraljice Marije 16, 11120 Belgrade, Serbia*

Correspondence should be addressed to Erdal Karapınar; erdalkarapinar@yahoo.com

We introduce the notion of θ-metric as a generalization of a metric by replacing the triangle inequality with a more generalized inequality. We investigate the topology of the spaces induced by a θ-metric and present some essential properties of it. Further, we give characterization of well-known fixed point theorems, such as the Banach and Caristi types in the context of such spaces.

## 1. Introduction

Celebrated Banach contraction mapping principle [1] can be considered as a revolution in fixed point theory and hence in nonlinear functional analysis. The statement of this well-known principle is simple, but the consequences are so strong: every contraction mapping in a complete metric space has a unique fixed point. Since fixed point theory has been applied to different sciences and also distinct branches of mathematics, this pioneer result of Banach has been generalized, extended, and improved in various ways in several abstract spaces. After that, many authors stated many types, generalizations, and applications of fixed point theory until now (see also [2, 3]).

In 1976, Caristi [4] defined an order relation in a metric space by using a functional under certain conditions and proved a fixed point theorem for such an ordered metric space. Denote $\mathbb{R}$ and $\mathbb{N}$ as the sets of all real and natural numbers, respectively.

The order relation is defined as follows.

**Lemma 1.** *Let $(X, d)$ be a metric space and $\phi : X \to \mathbb{R}$ a functional. Define the relation "$\preceq$" on $X$ by*

$$x \preceq y \iff d(x, y) \leq \phi(x) - \phi(y). \tag{1}$$

*Then, "$\preceq$" is a partial order relation on $X$ introduced by $\phi$, and $(X, d)$ is called an ordered metric space introduced by $\phi$. Apparently, if $x \preceq y$, then $\phi(x) \geq \phi(y)$.*

Caristi's fixed point theorem states that a mapping $T : X \to X$ has a fixed point provided that $(X, d)$ is a complete metric space and there exists a lower semicontinuous map $\phi : X \to \mathbb{R}$ such that

$$d(x, Tx) \leq \phi(x) - \phi(Tx), \quad \text{for every } x \in X. \tag{2}$$

This general fixed point theorem has found many applications in nonlinear analysis.

Many authors generalized Caristi's fixed point theorem and stated many types of it in complete metric spaces (see [5–8]). In particular, in 2010, Amini-Harandi [6] extended Caristi's fixed point and Takahashi's minimization theorems in complete metric space via the extension of partial ordered relation which is introduced in Lemma 1 and introduced some applications of such results.

One of the interesting generalizations of the notion of a metric is the concept of a fuzzy metric, given by Kramosil and Michálek [9], and Grabiec [10], independently. Later, George and Veeramani [11] investigated fuzzy metric structure and observed some important topological properties of such spaces. Furthermore, the authors [9–11] announced existence and uniqueness of a fixed point of certain mappings in the framework of such spaces.

*Definition 2* (see [11]). A binary operation $* : [0,1] \times [0,1] \to [0,1]$ is a continuous $t$-norm if it is a topological monoid with unit 1 such that $a * b \leq c * d$ whenever $a \leq c$ and $b \leq d$ ($a, b, c, d \in [0,1]$).

*Definition 3* (see [9, 10]). The 3-tuple $(X, M, *)$ is said to be a fuzzy metric space if $X$ is an arbitrary set, $*$ is a continuous $t$-norm, and $M : X \times X \times [0, +\infty) \to [0,1]$ is a fuzzy set satisfying the following conditions:

(1) $M(x, y, 0) = 0$,

(2) $M(x, y, t) = 1$ for all $t > 0$ if and only if $x = y$,

(3) $M(x, y, t) = M(y, x, t)$,

(4) $M(x, y, t) * M(y, z, t) \leq M(x, z, t+s)$,

(5) $M(x, y, \cdot) : [0, +\infty) \to [0, 1]$ is left continuous,

where $x, y, z \in X$ and $t, s > 0$.

In this paper, inspired from the definition of fuzzy metric spaces, we will introduce $\theta$-metric as an extension of metric spaces which is obtained by replacing the triangle inequality with a more generalized inequality. We also investigate the topology of the $\theta$-metric space and observe some fundamental properties of it. Furthermore, we give the characterization of the Banach and Caristi type fixed point theorems in the context of $\theta$-metric space.

First, we give the following definition.

*Definition 4.* Let $\theta : [0, +\infty) \times [0, +\infty) \to [0, +\infty)$ be a continuous mapping with respect to each variable. Let $\text{Im}(\theta) = \{\theta(s,t) : s \geq 0, t \geq 0\}$. A mapping $\theta$ is called an $B$-action if and only if it satisfies the following conditions:

(I) $\theta(0,0) = 0$ and $\theta(t,s) = \theta(s,t)$ for all $t, s \geq 0$,

(II)
$$\theta(s,t) < \theta(u,v) \quad \text{if } \begin{cases} \text{either } s < u, \ t \leq v, \\ \text{or } s \leq u, \ t < v, \end{cases} \quad (3)$$

(III) for each $r \in \text{Im}(\theta)$ and for each $s \in [0, r]$, there exists $t \in [0, r]$ such that $\theta(t, s) = r$,

(IV) $\theta(s, 0) \leq s$, for all $s > 0$.

We denote by $\Upsilon$ the set of all $B$-actions.

*Example 5.* The following functions are examples of $B$-action:

($\theta_1$) $\theta(t,s) = k(t+s)$, where $k \in (0, 1]$,

($\theta_2$) $\theta(t,s) = k(t+s+ts)$, where $k \in (0,1]$,

($\theta_3$) $\theta(t,s) = ts/(1+ts)$,

($\theta_4$) $\theta(s,t) = \sqrt{s^2 + t^2}$,

($\theta_5$) $\theta(t,s) = t + s + ts$,

($\theta_6$) $\theta(t,s) = t + s + \sqrt{ts}$,

($\theta_7$) $\theta(t,s) = (t+s)(1+ts)$.

Example 5 shows that the category of $B$-actions is uncountable. Next, we derive some lemmas which play a crucial role in our main results.

**Lemma 6.** *Let*

$$\Psi = \left\{ f : [0, +\infty) \longrightarrow [0, +\infty) \ : \ \begin{array}{c} f \text{ is continuous,} \\ \text{strictly increasing} \\ f(0) = 0 \\ f(t) < t \quad \forall t > 0. \end{array} \right\}. \quad (4)$$

*Then, there exists a correspondence between $\Upsilon$ and $\Psi$. In other words, $\Upsilon$ is an infinite set.*

*Proof.* For each $(t,s) \in [0, +\infty) \times [0, +\infty)$, define $\theta_f(t,s) = \lambda f(t+s)$, where $\lambda \in [0,1)$ and $f \in \Psi$. It is previous that $\theta_f \in \Upsilon$. Now, define

$$H : \Psi \longrightarrow \Upsilon$$
$$H(f) = \theta_f. \quad (5)$$

$H$ is well-defined and injective function which completes the proof. □

**Lemma 7.** *Let $\theta$ be a $B$-action. For each $r \in \text{Im}(\theta)$ and $s \in B = [0, r]$, there exist $t \in [0, r]$ and a function $\eta : [0, +\infty) \times [0, +\infty) \to [0, +\infty)$ such that $\eta(r, s) = t$. Then, one derives the following.*

($a_1$) $\eta(0, 0) = 0$.

($a_2$) $\theta(\eta(r, s), s) = r$ and $\theta(r, \eta(s, r)) = s$.

($a_3$) $\eta$ is continuous with respect to the first variable.

($a_4$) If $\eta(r, s) \geq 0$, then $0 \leq s \leq r$.

*Proof.* By (III) of Definition 4, for each $r \in \text{Im}(\theta)$ and for each $s \in [0, r]$, there exists $t \geq 0$ such that $\theta(t, s) = r$. Now, define $\eta(r, s) = t$. Let $t$ and $t'$ be two values such that $\eta(r, s) = t$ and $\eta(r, s) = t'$. If $t \neq t'$, then $t < t'$ or $t > t'$. If $t < t'$, then $r = \theta(t, s) < \theta(t', s) = r$ and this is a contradiction. For $t > t'$, we have the same argument. Thus, $\eta$ is well defined.

On the other hand, ($a_1$), ($a_2$), and ($a_4$) are straightforward from (III) of Definition 4. Also, ($a_3$) holds since if $r \in \text{Im}(\theta)$ and $\{r_n\}$ is a sequence in $\text{Im}(\theta)$ such that $r_n \to r$, then $\theta(\eta(r_n, \cdot), \cdot) = r_n$. Thus,

$$\theta\left(\lim_{n \to \infty} \eta(r_n, \cdot), \cdot\right) = \lim_{n \to \infty} \theta(\eta(r_n, \cdot), \cdot)$$
$$= r = \theta(\eta(r, \cdot), \cdot). \quad (6)$$

If $\lim_{n \to \infty} \eta(r_n, \cdot) > \eta(r, \cdot)$ or $\lim_{n \to \infty} \eta(r_n, \cdot) < \eta(r, \cdot)$, then by (II) of Definition 4 and (6), we conclude a contradiction. So, $\eta$ is continuous with respect to the first variable. □

The following definition arises from Lemma 7.

*Definition 8.* The function $\eta$, mentioned in Lemma 7, is called *B*-inverse action of $\theta$. One says that $\theta$ is regular if $\eta$ satisfies $\eta(r,r) = 0$, for each $r > 0$.

The set of all regular *B*-inverse actions will be denoted by $\Upsilon_r$.

*Example 9.* Let $\theta_1(t,s) = t+s$ and $\theta_2(t,s) = \sqrt[n]{t^n + s^n}$. It is evident that $\eta_1(t,s) = t-s$ and $\eta(t,s) = \sqrt[n]{t^n - s^n}$. Furthermore, $\eta_1$, $\eta_2$ satisfy all the conditions of Lemma 7. Note that $\theta_1$, $\theta_2$ are regular.

**Lemma 10.** *If $a \in X$, $c \in \operatorname{Im}(\theta)$ and $b \in [0,c]$, then $\theta(a,b) \le c$ implies that $a \le \eta(c,b)$.*

*Proof.* Suppose, on the contrary, that $\eta(c,b) < a$. Then, we have

$$c = \theta(\eta(c,b),b) < \theta(a,b), \tag{7}$$

a contradiction. □

## 2. Main Results

In this section, we introduce the notion $\theta$-metric and discuss the induced topology generated by $d_\theta$ on a nonempty subset $X$. In particular, we will show that $\tau_{d_\theta}$ is Hausdorff and first countable. Further, we derive that $(X, d_\theta)$ is a metrizable topological space.

*Definition 11.* Let $X$ be a nonempty set. A mapping $d_\theta : X \times X \to [0, +\infty)$ is called a $\theta$-metric on $X$ with respect to *B*-action $\theta \in \Upsilon$ if $d_\theta$ satisfies the following:

(A1) $d_\theta(x,y) = 0$ if and only if $x = y$,

(A2) $d_\theta(x,y) = d_\theta(y,x)$, for all $x, y \in X$,

(A3) $d_\theta(x,y) \le \theta(d_\theta(x,z), d_\theta(z,y))$, for all $x, y, z \in X$.

A pair $(X, d_\theta)$ is called a $\theta$-metric space.

*Remark 12.* If $(X, d_\theta)$ is a $\theta$-metric space, $\theta(s,t) = k(s+t)$, and $k \in (0,1]$, then $(X, d_\theta)$ is a metric space.

Conversely, for $\theta(s,t) = k(s+t)$, $k \in (0,1)$, we have that there exists metric space $(X, d)$ which is not $\theta$-metric space. For example, if $X = \{1, 2, 3\}$ and $d : X \times X \to [0, +\infty)$ defined by

$$d(1,2) = 1, \quad d(1,3) = 1, \quad d(2,3) = 2,$$

$$d(1,1) = 0, \quad d(2,2) = 0, \quad d(3,3) = 0, \quad k = \frac{1}{2}, \tag{8}$$

$d$ is a metric, but $d$ is not a $\theta$ metric.

*Example 13.* Let $X = \{x, y, z\}$ and $d_\theta : X \times X \to [0, +\infty)$ defined by

$$d_\theta(x,y) = 2, \quad d_\theta(x,z) = 6, \quad d_\theta(z,y) = 10,$$

$$d_\theta(x,x) = 0, \quad d_\theta(y,y) = 0, \quad d_\theta(z,z) = 0, \tag{9}$$

$$d_\theta(x,y) = d_\theta(y,x), \quad d_\theta(x,z) = d_\theta(z,x),$$

$$d_\theta(z,y) = d_\theta(z,x).$$

For $\theta(s,t) = s + t + st$, the function $d_\theta$ forms a metric, and hence, the pair $(X, d_\theta)$ is a $\theta$-metric space.

*Remark 14.* Notice that $d_\theta$, in Example 13, is not a metric on $X$, since we have $d_\theta(z,y) > d_\theta(x,y) + d_\theta(x,z)$.

*Definition 15.* Let $(X, d_\theta)$ be a $\theta$-metric space. An open ball $B_{d_\theta}(x,r)$ at a center $x \in X$ with a radius $r \in \operatorname{Im}(\theta)$ is defined as follows:

$$B_{d_\theta}(x,r) = \{y \in X : d_\theta(x,y) < r\}. \tag{10}$$

**Lemma 16.** *Every open ball is an open set.*

*Proof.* We show that, for each $x \in X$ and $r > 0$ and for each $y \in B_{d_\theta}(x,r)$, there exists $\delta > 0$ such that

$$B_{d_\theta}(y, \delta) \subset B_{d_\theta}(x,r). \tag{11}$$

By (III) of Definition 4, we can choose $\delta > 0$ such that $\theta(\delta, d_\theta(x,y)) = r$. Now, if $z \in B_{d_\theta}(y, \delta)$, then we have

$$d_\theta(z,x) \le \theta(d_\theta(z,y), d_\theta(y,x)) < \theta(\delta, d_\theta(y,x)) = r. \tag{12}$$

It means that $z \in B_{d_\theta}(x,r)$ and (11) is proved. □

**Lemma 17.** *If $(X, d_\theta)$ is a $\theta$-metric space, then the collection of open sets forms a topology, denoted by $\tau_{d_\theta}$. A pair $(X, \tau_\theta)$ is called topological space induced by a $\theta$-metric.*

**Lemma 18.** *The set $\{B_{d_\theta}(x, 1/n) : n \in \mathbb{N}\}$ is a local base at $x$, and the above topology is first countable.*

*Proof.* For each $x \in X$ and $r > 0$, we can find $n_0 \in \mathbb{N}$ such that $1/n_0 < r$. Thus, $B_{d_\theta}(x, 1/n) \subset B_{d_\theta}(x,r)$. This means that $\{B_{d_\theta}(x, 1/n) : n \in \mathbb{N}\}$ is a local base at $x$ and the above topology is first countable. □

**Theorem 19.** *A topological space $(X, \tau_{d_\theta})$ is Hausdorff.*

*Proof.* Let $x$, $y$ be two distinct points of $X$. Suppose that $0 < \alpha < d_\theta(x,y)$ is arbitrary. By Definition 4, we conclude that $\alpha \in \operatorname{Im}(\theta)$. Therefore, there exist $r, s > 0$ such that $\theta(r,s) = \alpha$. It is clear that $B_{d_\theta}(x,r) \cap B_{d_\theta}(x,s) = \emptyset$. For if there exists $z \in B_{d_\theta}(x,r) \cap B_{d_\theta}(x,s)$, then

$$d_\theta(x,y) \le \theta(d_\theta(x,z), d_\theta(z,y))$$
$$< \theta(r,s) = \alpha < d_\theta(x,y), \tag{13}$$

a contradiction. □

**Theorem 20.** *Let $(X, d_\theta)$ be a $\theta$-metric space and $\tau_{d_\theta}$ the topology induced by the $\theta$-metric. Then, for a sequence $\{x_n\}$ in $X$, $x_n \to x$ if and only if $d_\theta(x_n, x) \to 0$ as $n \to \infty$.*

*Proof.* Suppose that $x_n \to x$. Then, for each $\epsilon > 0$, there exists $n_0 \in \mathbb{N}$ such that $x_n \in B_{d_\theta}(x, \epsilon)$, for all $n \geq n_0$. Thus, $d_\theta(x_n, x) < \epsilon$; that is, $d_\theta(x_n, x) \to 0$ as $n \to \infty$. The converse is verified easily. □

**Theorem 21.** *Let $(X, d_\theta)$ be a $\theta$-metric space and $x_n \to x$, $y_n \to y$. Then,*

$$d_\theta(x_n, y_n) \longrightarrow d_\theta(x, y). \tag{14}$$

*Proof.* For each $n \in \mathbb{N}$, there exists $K > 0$ such that, for all $n \geq K$,

$$d_\theta(x_n, x) < \frac{1}{n}, \qquad d_\theta(y_n, y) < \frac{1}{n}. \tag{15}$$

Thus, by the continuity of $\theta$ with respect to each variable, we have

$$d_\theta(x, y) \leq \theta(d_\theta(x, x_n), \theta(d_\theta(x_n, y_n), d_\theta(y_n, y))) \\ < \theta\left(\frac{1}{n}, \theta\left(d_\theta(x_n, y_n), \frac{1}{n}\right)\right). \tag{16}$$

Therefore,

$$\begin{aligned} d_\theta(x, y) &\leq \lim_{n \to \infty} \theta\left(\frac{1}{n}, \theta\left(d_\theta(x_n, y_n), \frac{1}{n}\right)\right) \\ &= \theta\left(0, \theta\left(\lim_{n \to \infty} d_\theta(x_n, y_n), 0\right)\right) \\ &\leq \lim_{n \to \infty} d_\theta(x_n, y_n) \\ &\leq \lim_{n \to \infty} \theta(d_\theta(x, x_n), \theta(d_\theta(x, y), d_\theta(y_n, y))) \\ &\leq \lim_{n \to \infty} \theta\left(\frac{1}{n}, \theta\left(d_\theta(x, y), \frac{1}{n}\right)\right) \\ &\leq d_\theta(x, y). \end{aligned} \tag{17}$$

Thus, $d_\theta(x_n, y_n) \to d_\theta(x, y)$. □

**Lemma 22.** *Let $(X, d_\theta)$ be a $\theta$-metric space. Let $\{x_n\}$ be a sequence in $X$ and $x_n \to x$. Then, $x$ is unique.*

*Proof.* Suppose that $x_n \to x$ and $x_n \to y$. We show that $x = y$. For each $n \in \mathbb{N}$, there exists $N > 0$ such that $d_\theta(x_n, x) < 1/n$ and $d_\theta(x_n, y) < 1/n$. By the continuity of $\theta$, we have

$$0 \leq d_\theta(x, y) \leq \theta(d_\theta(x_n, x), d_\theta(x_n, y)) \\ < \theta\left(\frac{1}{n}, \frac{1}{n}\right) \longrightarrow 0, \quad (n \longrightarrow \infty). \tag{18}$$

Hence, we have $x = y$. □

**Definition 23.** Let $(X, d_\theta)$ be a $\theta$-metric space. Then, for a sequence $\{x_n\}$ in $X$, one says that $\{x_n\}$ is a Cauchy sequence if for each $\epsilon > 0$, there exists $N > 0$ such that, for all $m \geq n \geq N$, $d_\theta(x_n, x_m) < \epsilon$.

**Definition 24.** Let $(X, d_\theta)$ be a $\theta$-metric space. One says that $(X, d_\theta)$ is complete $\theta$-metric space if every Cauchy sequence $\{x_n\}$ is convergent in $X$.

**Lemma 25** (see [12]). *A Hausdorff topological space $(X, \tau_{d_\theta})$ is metrizable if and only if it admits a compatible uniformity with a countable base.*

In the following theorems we apply the previous lemma and the concept of uniformity (see [12] for more information) to prove the metrizability of a topological space $(X, \tau_{d_\theta})$.

**Theorem 26.** *Let $(X, d_\theta)$ be a $\theta$-metric space. Then, $(X, \tau_{d_\theta})$ is a metrizable topological space.*

*Proof.* For each $n \in \mathbb{N}$, define

$$\mathcal{U}_n = \left\{(x, y) \in X \times X : d_\theta(x, y) < \frac{1}{n}\right\}. \tag{19}$$

We will prove that $\{\mathcal{U}_n : n \in \mathbb{N}\}$ is a base for a uniformity $\mathcal{U}$ on $X$ whose induced topology coincides with $\tau_{d_\theta}$.

We first note that for each $n \in \mathbb{N}$,

$$\{(x, x) : x \in X\} \subseteq \mathcal{U}_n, \quad \mathcal{U}_{n+1} \subseteq \mathcal{U}_n, \quad \mathcal{U}_n = \mathcal{U}_n^{-1}. \tag{20}$$

On the other hand, for each $n \in \mathbb{N}$, there is, by the continuity of $\theta$, an $m \in \mathbb{N}$ such that

$$m > 2n, \qquad \theta\left(\frac{1}{m}, \frac{1}{m}\right) < \frac{1}{n}. \tag{21}$$

Then, $\mathcal{U}_m \circ \mathcal{U}_m \subseteq \mathcal{U}_n$: indeed, let $(x, y) \in \mathcal{U}_m$ and $(y, z) \in \mathcal{U}_m$. Thus,

$$d_\theta(x, z) \leq \theta(d_\theta(x, y), d_\theta(y, z)) \\ < \theta\left(\frac{1}{m}, \frac{1}{m}\right) < \frac{1}{n}. \tag{22}$$

Therefore, $(x, z) \in \mathcal{U}_n$. Hence, $\{\mathcal{U}_n : n \in \mathbb{N}\}$ is a base for a uniformity $\mathcal{U}$ on $X$. Since for each $x \in X$ and each $n \in \mathbb{N}$, $\mathcal{U}_n(x) = \{y \in X : d_\theta(x, y) < 1/n\}$. We deduce from Lemma 25 that $(X, \tau_{d_\theta})$ is a metrizable topological space. □

Let us recall that a metrizable topological space $(X, \tau)$ is said to be completely metrizable if it admits a complete metric [13].

**Theorem 27.** *Let $(X, d_\theta)$ be a complete $\theta$-metric space. Then, $(X, \tau_{d_\theta})$ is completely metrizable.*

*Proof.* It follows from the proof of Theorem 26 that $\{\mathcal{U}_n : n \in \mathbb{N}\}$ is a base for a uniformity $\mathcal{U}$ on $X$ compatible with $\tau_{d_\theta}$, where $\mathcal{U}_n = \{(x, y) \in X \times X : d_\theta(x, y) < 1/n\}$, for every $n \in \mathbb{N}$. Then, there exists a metric $d$ on $X$ whose induced uniformity coincides with $\mathcal{U}$. We want to show that the metric $d$ is complete. Indeed, given a Cauchy sequence $\{x_n\}$ in $(X, d)$, we will prove that $\{x_n\}$ is a Cauchy sequence in $(X, d)$. To this end, fix $\epsilon > 0$. Choose $k \in \mathbb{N}$ such that $1/k < \epsilon$. Then, there exists $n_0 \in \mathbb{N}$ such that $(x_m, x_n) \in \mathcal{U}_k$ for every $n, m \geq n_0$. Consequently, for each $n, m \geq n_0$, $d_\theta(x_n, x_m) \leq 1/k < \epsilon$. We have shown that $\{x_n\}$ is a Cauchy sequence in the complete $\theta$-metric space $(X, d)$ and so is convergent with respect to $(X, d)$. Thus, $(X, d)$ is a complete metric space. □

## 3. Two Fixed Point Theorems

In this section, we introduce two fixed point theorems in $\theta$-metric spaces. First, we introduce the Banach fixed point and Caristi's fixed point theorems in such spaces.

### 3.1. Banach Fixed Point Theorem

**Theorem 28.** *Let $(X, d_\theta)$ be a complete $\theta$-metric space and $f : X \to X$ a mapping that satisfies the following:*

$$d_\theta(fx, fy) \leq \alpha d_\theta(x, y), \tag{23}$$

*for each $x, y \in X$, where $\alpha \in [0, 1)$. Then, $f$ has a unique fixed point.*

*Proof.* Let $x_0 \in X$ and $x_{n+1} = fx_n$. We divide our proof into four steps.

*Step 1.* We claim that $d_\theta(x_n, x_{n+1}) \to 0$. Indeed, we have

$$\begin{aligned} d_\theta(x_{n+1}, x_n) &\leq \alpha d_\theta(x_n, x_{n-1}) \\ &\leq \alpha^2 d_\theta(x_{n-1}, x_{n-2}) \\ &\vdots \\ &\leq \alpha^n d_\theta(x_1, x_0) \to 0, \quad \text{as } (n \to \infty). \end{aligned} \tag{24}$$

Thus, we have $d_\theta(x_n, x_{n+1}) \to 0$.

*Step 2.* We assert that the sequence $\{x_n\}$ is bounded. Suppose, on the contrary, that $\{x_n\}$ is an unbounded sequence. Thus, there exists subsequence $\{n(k)\}$ such that $n(1) = 1$ and for each $k \in \mathbb{N}$, $n(k+1)$ is minimal in the sense that the relation

$$d_\theta(x_{n(k+1)}, x_{n(k)}) > 1 \tag{25}$$

does not hold and

$$d_\theta(x_m, x_{n(k)}) \leq 1 \tag{26}$$

holds for all $m \in \{n(k)+1, n(k)+2, \ldots, n(k+1)-1\}$. Hence, by using the triangle inequality, we derive that

$$\begin{aligned} 1 &< d_\theta(x_{n(k+1)}, x_{n(k)}) \\ &\leq \theta(d_\theta(x_{n(k+1)}, x_{n(k+1)-1}), d_\theta(x_{n(k+1)-1}, x_{n(k)})) \\ &\leq \theta(d_\theta(x_{n(k+1)}, x_{n(k+1)-1}, 1)). \end{aligned} \tag{27}$$

By taking the limit from two sides of (27) and using (II) of Definition 4, we derive that

$$d_\theta(x_{n(k+1)}, x_{n(k)}) \to 1^+ \quad \text{as } (k \to +\infty). \tag{28}$$

Also, we have

$$\begin{aligned} 1 &< d_\theta(x_{n(k+1)}, x_{n(k)}) \\ &\leq d_\theta(x_{n(k+1)-1}, x_{n(k)-1}) \\ &\leq \theta(d_\theta(x_{n(k+1)-1}, x_{n(k)}), d_\theta(x_{n(k)}, x_{n(k)-1})) \\ &\leq \theta(1, d_\theta(x_{n(k)}, x_{n(k)-1})), \end{aligned} \tag{29}$$

which implies that

$$d_\theta(x_{n(k+1)-1}, x_{n(k)-1}) \to 1^+ \quad \text{as } (k \to +\infty). \tag{30}$$

Since $d_\theta(x_{n(k+1)}, x_{n(k)}) \leq \alpha d_\theta(x_{n(k+1)-1}, x_{n(k)-1})$, we have $1 \leq \alpha 1$, a contradiction. Thus, the sequence $\{x_n\}$ is bounded.

*Step 3.* We will show that $\{x_n\}$ is a Cauchy sequence. Let $m, n \in \mathbb{N}$ with $m > n$

$$\begin{aligned} d_\theta(x_m, x_n) &= d_\theta(fx_{m-1}, fx_{n-1}) \\ &\leq \alpha d_\theta(fx_{m-2}, fx_{n-2}) \\ &\vdots \\ &\leq \alpha^n d_\theta(fx_{m-n}, x_0). \end{aligned} \tag{31}$$

Since $\{x_n\}$ is a bounded sequence, therefore, $\lim_{n,m \to \infty} d_\theta(x_m, x_n) = 0$; that is, $\{x_n\}$ is a Cauchy sequence. Thus, there exists $x \in X$ such that $x_n \to x$. Further, we derive that

$$\begin{aligned} d_\theta(x_{n+1}, fx) &= d_\theta(fx_n, fx) \\ &\leq \alpha d_\theta(x_n, x)) \to 0, \quad (n \to \infty). \end{aligned} \tag{32}$$

It means that $x_{n+1} \to fx$; that is, $fx = x$.

*Step 4.* In the last step, we prove that the $x$ is the unique fixed point of $T$. Suppose, on the contrary, that $x, y$ are two distinct fixed points of $f$. So, we get that

$$d_\theta(y, x) = d_\theta(fy, fx) \leq \alpha d_\theta(y, x) \tag{33}$$

is a contradiction. This completes the proof. □

### 3.2. Caristi-Type Fixed Point Theorem

*Definition 29.* Suppose that $(X, d_\theta)$ be a complete $\theta$-metric space and $\mathscr{P}$ the class of all maps $\psi : X \times X \to [0, +\infty)$ which satisfies the following conditions:

($E_1$) there exists $\hat{x} \in X$ such that $\psi(\hat{x}, \cdot)$ is bounded below and lower semicontinuous, and $\psi(\cdot, y)$ is upper semicontinuous for each $y \in X$,

($E_2$) $\psi(x, x) = 0$, for each $x \in X$,

($E_3$) $\theta(\psi(x, y), \psi(y, z)) \leq \psi(x, z)$, for each $x, y, z \in X$.

**Lemma 30.** *By Definition 29, one has*

$$\psi(x, y) \leq \eta(\psi(x, z), \psi(y, z)), \tag{34}$$

*for each $x, y, z \in X$.*

*Proof.* By Lemma 10, we obtain the desired result. □

*Example 31.* Let $\theta(t, s) = ts/(1+ts)$; thus, $\text{Im}(\theta) = [0, 1)$. Now, let $\phi : X \to \mathbb{R}$ be a lower bounded, lower semicontinuous function and

$$\psi(x, y) = \begin{cases} e^{\phi(y) - \phi(x)} & x \neq y \\ 0 & x = y. \end{cases} \tag{35}$$

Then, $\psi$ satisfies all conditions of Definition 29.

*Example 32.* Let $\theta(t,s) = \sqrt[2n+1]{t^{2n+1} + s^{2n+1}}$; thus, $\text{Im}(\theta) = [0, +\infty)$. Now, let $\phi : X \to \mathbb{R}$ be a lower bounded, lower semicontinuous function and

$$\psi(x, y) = \sqrt[2n+1]{\phi(y) - \phi(x)}. \tag{36}$$

Then, $\psi$ satisfies all conditions of Definition 29. Also, $\eta(t,s) = \sqrt[2n+1]{t^{2n+1} - s^{2n+1}}$, and $\theta$ is regular.

From now to end, we assume that $\theta$ is regular (see Definition 8).

**Lemma 33.** *Let $(X, d_\theta)$ be a complete $\theta$-metric space and $\psi \in \mathcal{P}$. Let $\gamma : [0, +\infty) \to [0, +\infty)$ be $\theta$-subadditive; that is, $\gamma(\theta(x, y)) \le \theta(\gamma(x), \gamma(y))$, for each $x, y \in [0, +\infty)$, nondecreasing continuous map such that $\gamma^{-1}\{0\} = \{0\}$. Define the order $\prec$ on $X$ by*

$$x \prec y \iff \gamma(d_\theta(x, y)) \le \psi(x, y), \tag{37}$$

*for any $x, y \in X$. Then, $(X, \prec)$ is a partial order set which has minimal elements.*

*Proof.* At first, we show that $(X, \prec)$ is a partial ordered set. For each $x \in X$, we have $0 = \gamma(0) = \gamma(d_\theta(x, x)) \le \psi(x, x) = 0$. Thus, $x \prec x$. If $x \prec y$ and $y \prec x$, then $\gamma(d_\theta(x, y)) \le \psi(x, y)$ and $\gamma(d_\theta(y, x)) \le \psi(y, x)$. Thus, we give

$$\gamma(\theta(d_\theta(x, y), d_\theta(x, y))) \le \theta(\gamma(d_\theta(x, y)), \gamma(d_\theta(x, y)))$$
$$\le \theta(\psi(x, y), \psi(y, x))$$
$$\le \psi(x, x) = 0. \tag{38}$$

It means that $x = y$. Finally, if $x \prec y$ and $y \prec z$, then $\gamma(d_\theta(x, y)) \le \psi(x, y)$ and $\gamma(d_\theta(y, z)) \le \psi(y, z)$. Thus, we give

$$\gamma(d_\theta(x, z)) \le \gamma(\theta(d_\theta(x, y), d_\theta(y, z)))$$
$$\le \theta(\gamma(d_\theta(x, y)), \gamma(d_\theta(y, z)))$$
$$\le \theta(\psi(x, y), \psi(y, z))$$
$$\le \psi(x, z). \tag{39}$$

It means that $x \prec z$. Thus, $(X, \prec)$ is a partial ordered set.

To show that $(X, \prec)$ has minimal elements, we show that any decreasing chain has a lower bound. Indeed, let $\{x_\alpha\}_{\alpha \in \Gamma}$ be a decreasing chain; then we have

$$0 \le \gamma(d_\theta(x_\alpha, x_\beta))$$
$$\le \psi(x_\alpha, x_\beta) \tag{40}$$
$$\le \eta(\psi(\hat{x}, x_\beta), \psi(\hat{x}, x_\alpha));$$

by definition of $\eta$, we have $\psi(\hat{x}, x_\alpha) \le \psi(\hat{x}, x_\beta)$. Thus, $\{\psi(\hat{x}, x_\alpha)\}_{\alpha \in \Gamma}$ is decreasing net of reals which is bounded below. Let $\{\alpha_n\}$ be an increasing sequence of elements from $\Gamma$ such that

$$\lim_{n \to \infty} \psi(\hat{x}, x_{\alpha_n}) = \inf\{\psi(\hat{x}, x_\alpha) : \alpha \in \Gamma\} = \rho. \tag{41}$$

Then, for each $m \ge n$, we infer that

$$\gamma(d_\theta(x_{\alpha_n}, x_{\alpha_m})) \le \psi(x_{\alpha_n}, x_{\alpha_m})$$
$$\le \eta(\psi(\hat{x}, x_{\alpha_n}), \psi(\hat{x}, x_{\alpha_m})). \tag{42}$$

By taking limit from two sides of (42), the regularity of $\theta$, and continuity of $\eta$, we give

$$\limsup_{n,m \to \infty} \gamma(d_\theta(x_{\alpha_n}, x_{\alpha_m}))$$
$$\le \psi(x_{\alpha_n}, x_{\alpha_m})$$
$$\le \limsup_{n,m \to \infty} \eta(\psi(\hat{x}, x_{\alpha_n}), \psi(\hat{x}, x_{\alpha_m})) \tag{43}$$
$$\le \eta(\rho, \rho) = 0.$$

Then, our assumption on $\gamma$ implies that $\{x_{\alpha_n}\}$ is a Cauchy sequence and therefore converges to some $x \in X$. Since $\gamma$ is continuous and $\psi(\cdot, x_{\alpha_n})$ is upper semicontinuous, then we have

$$\gamma(d_\theta(x, x_{\alpha_n})) = \limsup_{m \to \infty} \gamma(d_\theta(x_{\alpha_m}, x_{\alpha_n}))$$
$$\le \limsup_{m \to \infty} \psi(x_{\alpha_m}, x_{\alpha_n}) \tag{44}$$
$$\le \psi(x, x_{\alpha_n}).$$

This shows that $x \prec x_{\alpha_n}$ for all $n \ge 1$, which means that $x$ is lower bound for $\{x_{\alpha_n}\}$. In order to see that $x$ is also a lower bound for $\{x_\alpha\}_{\alpha \in \Gamma}$, let $\beta \in \Gamma$ be such that $x_\beta \prec x_{\alpha_n}$ for all $n \ge 1$. Then, for each $n \in \mathbb{N}$, we have

$$0 \le \gamma(d_\theta(x_\beta, x_{\alpha_n}))$$
$$\le \psi(x_\beta, x_{\alpha_n}) \tag{45}$$
$$\le \eta(\psi(\hat{x}, x_{\alpha_n}), \psi(\hat{x}, x_\beta)).$$

Hence, for all $n \ge 1$,

$$\psi(\hat{x}, x_\beta) \le \psi(\hat{x}, x_{\alpha_n}) \tag{46}$$

which implies that

$$\psi(\hat{x}, x_\beta) = \inf\{\psi(\hat{x}, x_\alpha) : \alpha \in \Gamma\}$$
$$= \lim_{n \to \infty} \psi(\hat{x}, x_{\alpha_n}). \tag{47}$$

Thus, from (46), we get $\lim_{n \to \infty} x_{\alpha_n} = x_\beta$, which implies that $x_\beta = x$. Therefore, for any $\alpha \in \Gamma$, there exists $n \in \mathbb{N}$ such that $x_{\alpha_n} \prec x_\alpha$; that is, $x$ is a lower bound of $\{x_\alpha\}$. Zorn's lemma will therefore imply that $(X, \prec)$ has minimal elements. □

**Theorem 34.** *Let $(X, d_\theta)$ be a complete $\theta$-metric space and $\psi \in \mathcal{P}$. Let $\gamma : [0, +\infty) \to [0, +\infty)$ be as in Lemma 33. Let $T : X \to X$ be a map satisfying the following:*

$$\gamma(d_\theta(x, Tx)) \le \psi(Tx, x), \tag{48}$$

*for any $x \in X$. Then, $T$ has a fixed point.*

*Proof.* By Lemma 33, $(X, \prec)$ has a minimal element say $\overline{x}$. Thus, $T\overline{x} \prec \overline{x}$. It means that $T\overline{x} = \overline{x}$. □

**Corollary 35.** *Let $(X, d_\theta)$ be a complete $\theta$-metric space and $\psi \in \mathcal{P}$. Let $\gamma : [0, +\infty) \to [0, +\infty)$ be as in Lemma 33. Let $T : X \to 2^X$ be a multivalued mapping satisfying the following:*

$$\gamma(d_\theta(x, y)) \le \psi(y, x), \quad \forall y \in Tx. \tag{49}$$

*Then, $T$ has an endpoint; that is, there exists $\overline{x} \in X$ such that $T\overline{x} = \{\overline{x}\}$.*

In Corollary 35, we can introduce many types of Caristi's fixed point theorem as follows.

If we set $\psi$ as in Example 32, then (48) has the following form:

$$\gamma(d_\theta(x, Tx)) \le \sqrt[2n+1]{\phi(Tx) - \phi(x)}. \tag{50}$$

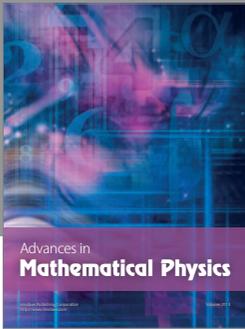 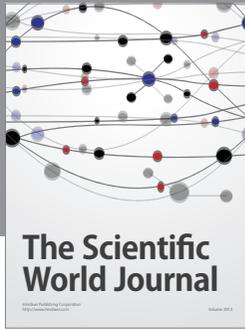 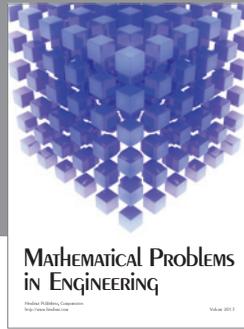 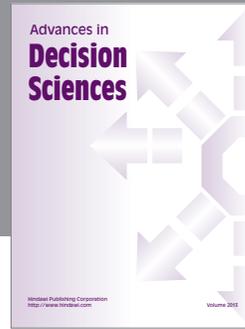 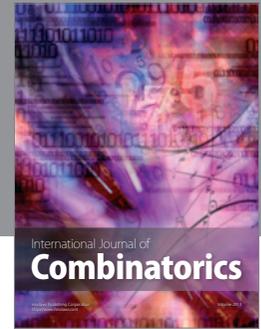
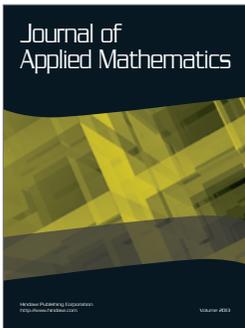 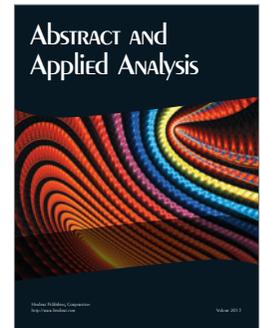
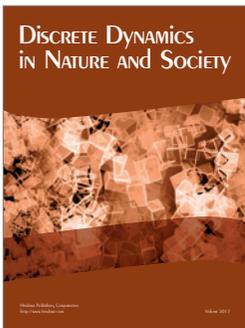 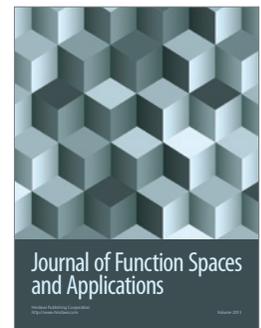
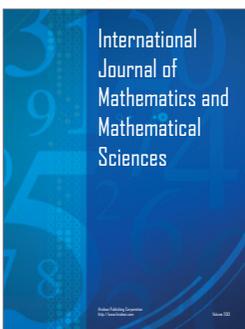 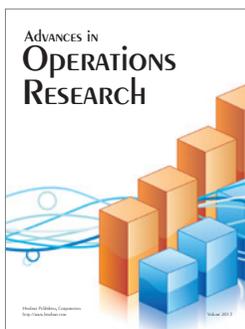 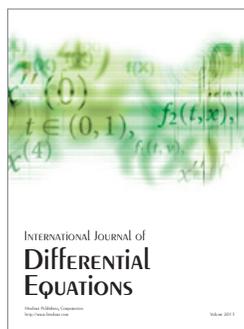 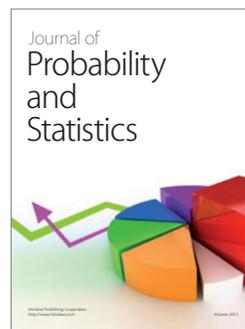 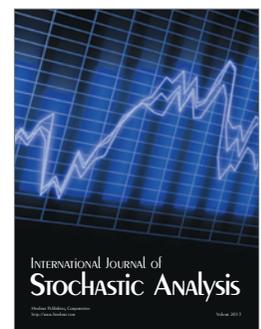
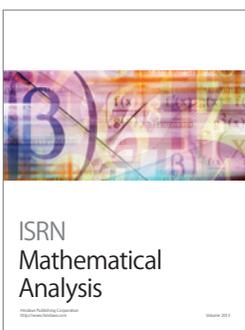 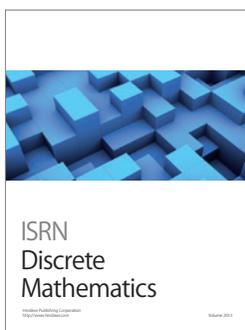 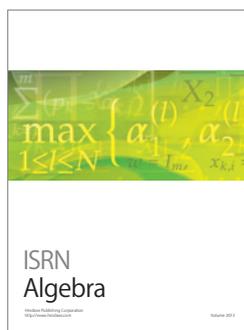 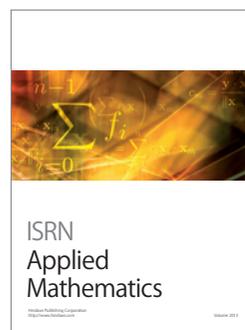 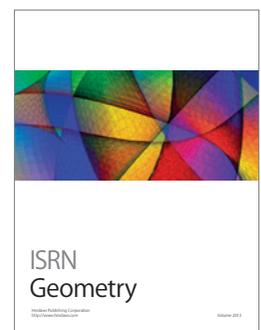

Submit your manuscripts at
http://www.hindawi.com